\documentclass[12pt]{article}
\usepackage[table]{xcolor} 
\definecolor{lightblue}{rgb}{0.7,0.8,1} 
\usepackage{latexsym}
\usepackage{amssymb}
\usepackage[noadjust]{cite}
\usepackage{tikz}
\usetikzlibrary{matrix,arrows}
\usepackage{tikz-cd}
\usepackage{circuitikz}
\usetikzlibrary{circuits}
\usetikzlibrary[circuits]
\usetikzlibrary{circuits.logic.US,circuits.logic.IEC}
\usepackage{graphicx}
\usepackage{wrapfig}
\usepackage{booktabs}
\usepackage{url} 
\usepackage{subcaption}
\usepackage{amsmath}
\usepackage{float}
\usepackage{parskip}
\usepackage{nohyperref} 
\usepackage{array}    
\renewcommand{\arraystretch}{1.5} 
\Large 

\usepackage [a4paper,left=2.5cm,right=2.5cm,top=2.5cm, bottom=2.5cm]{geometry}

\makeatletter
\def\section{\@startsection {section}{1}{\z@}{-3.5ex plus -1ex minus
 -.2ex}{2.3ex plus .2ex}{\large\bf}}
\def\subsection{\@startsection{subsection}{2}{\z@}{-3.25ex plus -1ex minus
 -.2ex}{1.5ex plus .2ex}{\normalsize\bf}}
\makeatother

\makeatletter
\@addtoreset{equation}{section}

\newcommand{\nc}{\newcommand}
\newcommand{\rnc}{\renewcommand}

\nc{\be}{\bea}
\nc{\ee}{\eea}

\rnc{\a}{\alpha}
\nc{\ab}{\bar{\a}}
\nc{\ap}{\a^{+}}
\nc{\abm}{\ab^{-}}
\rnc{\b}{\beta}
\nc{\bb}{\bar{\b}}
\nc{\bbp}{\bb_{\zb}^{+}}
\nc{\bm}{\b_{z}^{-}}
\nc{\oa}{\overline{\a}}
\nc{\ob}{\overline{\b}}
\rnc{\gg}{\gamma}
\rnc{\d}{\delta}
\nc{\f}{\phi}
\nc{\fb}{\bar{\phi}}
\nc{\vf}{\varphi}
\nc{\p}{\psi}

\rnc{\c}{\chi}
\nc{\la}{\lambda}
\nc{\m}{\mu}
\nc{\n}{\nu}
\rnc{\o}{\omega}
\nc{\Om}{\Omega}
\rnc{\t}{\theta}
\nc{\eps}{\epsilon}
\rnc{\S}{\Sigma}
\nc{\F}{\Phi}

\nc{\trac}[2]{{\textstyle\frac{#1}{#2}}}

\nc{\ex}[1]{\mbox{e}^{\,\textstyle#1}}

\nc{\mat}[4]{\left(\begin{array}{cc}#1&#2\\#3&#4\end{array}\right)}

\nc{\som}[9]{\left(\begin{array}{ccc}#1&#2&#3\\#4&#5&#6\\#7&#8&#9%
\end{array}\right)}

\nc{\tr}{\mathop{\mbox{tr}}\nolimits}
\nc{\ad}{\mathop{\mbox{ad}}\nolimits}
\nc{\Tr}{\mathop{\mbox{Tr}}\nolimits}
\nc{\Det}{\mathop{\mbox{Det}}\nolimits}
\nc{\rk}{\mathop{\mbox{rk}}\nolimits}
\nc{\ra}{\rightarrow}
\nc{\Ra}{\Rightarrow}
\nc{\LRa}{\Leftrightarrow}
\nc{\ot}{\otimes}
\rnc{\ss}{\subset}
\nc{\nul}{\noindent\underline}
\nc{\non}{\nonumber\\}

\nc{\subs}[1]{{\vspace*{0.5cm}}%
{\noindent\underline{#1}}{\addcontentsline{toc}{subsection}{#1}}%
{\vspace*{0.3cm}}}
\rnc{\subs}[1]{\subsection{#1}}

\usepackage{latexsym}
\usepackage{mathrsfs}
\usepackage{amsfonts}
\usepackage{amsbsy}
\usepackage{indentfirst}
\usepackage{multirow}
\usepackage[Glenn]{fncychap}
\begin{document}
\begin{titlepage}
\renewcommand{\thefootnote}{\fnsymbol{footnote}}
\begin{center}
{\large \bf Association Rules Machine Learning complete intersection Calabi-Yau $5$-Folds and $6$-Folds}
\vspace{1.0cm}

{Kaniba Mady Keita$^{a,b}$}\footnote{E-mail: madyfalaye@gmail.com, kanibamady.keita@usttb.edu.ml}

\vspace{.5cm}
{${a}$ }{\it\small  Centre de Calcul de Mod\'{e}lisation et de Simulation: CCMS\\
 Department of Physics, Faculty of Sciences and Techniques, University of Sciences, Techniques and Technologies of Bamako, FST-USTTB, BP:E3206, Mali.}\\
$b$ {\it\small Centre de Recherche en Phyique Quantique et de ses Applications: CRPQA, Bamako, Mali.}\\

\vspace{.3cm}
\today
\end{center}

\vspace{1.5cm}

\centerline{\textbf{Abstract}}\vspace{0.5cm}

Association rule machine learning is applied to the dataset of complete intersection Calabi--Yau 5-folds and 6-folds in order to uncover hidden patterns among their Hodge numbers. These Hodge numbers—six for the 5-folds and nine for the 6-folds—serve as the items in our analysis. For the 5-folds, we discover 60 significant association rules. For example, within the dataset, if $h^{1,3} = 0$ and $h^{2,2} = 5$, then $h^{1,1} = 3$ with $99.43\%$ confidence. Similarly, if $h^{2,1} = 0$, $h^{1,3} = 0$, and $h^{2,2} = 5$, then $h^{1,1} = 3$ with $99.42\%$ confidence. For the 6-folds, we identify 160 association rules across a dataset of $1,482,022$ examples. A particularly striking observation is that $h^{1,2} = h^{1,3} = h^{1,4} = h^{2,3} = 0$ for all entries in this dataset. These types of association rules are especially valuable because the Hodge numbers of complete intersection Calabi–Yau 5-folds have only been computed for approximately $53\%$ of the dataset, while those of 6-folds remain largely undetermined. The discovered patterns provide predictive insights that can guide future computations and theoretical developments.
\end{titlepage}
\setcounter{footnote}{0}

\section{Introduction}

Topological invariants of complete intersection Calabi–Yau three-folds (CICY3s) have been computed, and the corresponding dataset consists of $7890$ CICY3s~\cite{pgreen}. Similarly, there are $921{,}497$ known complete intersection Calabi–Yau four-folds~\cite{gray}. In the case of complete intersection Calabi–Yau five-folds (CICY5s), $27{,}068$ spaces were obtained by the authors of~\cite{alawadhi}, but the cohomological data has been computed for only $12{,}433$ of them (approximately $53.7\%$). 

The space of Calabi–Yau six-folds is even more intriguing. Edward Hirst and Tancredi Schettini Gherardini identified approximate Hodge numbers for a dataset containing $1{,}482{,}022$ Calabi–Yau six-fold candidates with a total weight sum less than $200$~\cite{Edward}. These Hodge numbers were shown to match exactly with the values computed by V.~N.~Dumachev~\cite{dumachev}. Additionally, these six-dimensional weighted projective spaces with transverse polynomials were confirmed to be consistent by Maximilian Kreuzer and Harald Skarke~\cite{Skarke}.  

In recent years, machine learning has been successfully applied to the study of complete intersection Calabi–Yau manifolds~\cite{kaniba, kaniba2, Yang-hui, yang, yang2, wei, harold, Ashmore, daata, berman, rehan, lukas, rue, bull1, bull2, erbin, erbin2}. These approaches are believed to simplify the computation of topological invariants of such manifolds. Most of the machine learning applications in this context have focused on classification, regression, and clustering algorithms. While these techniques are quite effective for fully determined datasets, this is unfortunately not the case for Calabi–Yau $5$-folds and $6$-folds.

For incomplete datasets, one particularly useful technique is association rule machine learning (also known as association rule data mining)~\cite{Zhangass}, which can be employed to propose logical rules for estimating the missing data entries.

In this work, we apply association rule machine learning using the Apriori algorithm~\cite{Srikant} and some of its improved versions (see, for instance,~\cite{Juanxia, Maolegi, Zawaidah}). Our implementation uses the \texttt{arulesViz} package~\cite{Hahsler} in RStudio. The main evaluation metrics for an association rule are the \emph{support}, \emph{confidence}, and \emph{lift} of the rule. These metrics will be defined in Section 2.

The remainder of this paper is organized as follows. In Section~\ref{sec:background}, we introduce the basic theoretical concepts of association rules and describe how they are implemented using the \texttt{arulesViz} package. Section~\ref{sec:motivation} discusses the challenges involved in computing the cohomological datasets for Calabi–Yau $5$-folds and $6$-folds and motivates the application of association rule techniques. Section~\ref{sec:results} presents the results of applying the Apriori algorithm to these datasets. Finally, Section~\ref{sec:conclusion} offers concluding remarks and outlines directions for future work.

\section{Basic Theoretical Concepts of Association Rules}\label{sec:background}
The basic concepts we present here can be found in \cite{jiawei, aze}.

Let us denote our dataset by \( \mathscr{D} \), where each row corresponds to a transaction. The first column of \( \mathscr{D} \) represents the identity of each transaction, while the remaining columns correspond to items.

In our context, each transaction represents a specific Calabi–Yau manifold, and the identity refers to its sequential label within the dataset. For instance, the CICY3 dataset contains \( 7{,}890 \) such identities, whereas the CICY6 dataset contains \( 1{,}482{,}022 \).

The items in the dataset are the individual Hodge numbers \( h^{i,j} \) associated with each manifold. By definition, an \emph{itemset} \( X \) is a collection of one or more such items. If \( X \) contains \( k \) items, it is referred to as a \emph{\( k \)-itemset}.

A \textbf{association rule} is a logical implication between two itemsets, typically expressed as:
\[
X \Rightarrow Y,
\]
which means that if the itemset \( X \) appears in a transaction, then the itemset \( Y \) is likely to appear as well. Here, \( X \) is called the \textbf{left-hand side (LHS)} and \( Y \) is called the \textbf{right-hand side (RHS)} of the rule.

When \( X = \emptyset \), the association rule
\[
\emptyset \Rightarrow Y
\]
means that the itemset \( Y \) is likely to appear in all transactions, regardless of the presence of any other items. In other words, the rule holds unconditionally—\( Y \) is always true.

The \textbf{support} of the association rule \( X \Rightarrow Y \), denoted \( \text{support}(X \Rightarrow Y) \), is the fraction of transactions in the dataset \( \mathscr{D} \) that contain both itemsets \( X \) and \( Y \); that is, transactions containing \( X \cup Y \).

The \emph{confidence} of the rule, denoted by \( \mathrm{conf}(X \Rightarrow Y) \), is, by definition, the conditional probability \( p(Y \mid X) \). In other words, confidence is the proportion of transactions in \( \mathscr{D} \) that contain \( X \) and also contain \( Y \).
A few statistical parameters for association rules are given in Table~\ref{tab:assoc_metrics} below:
\begin{table}[h]
\centering
\begin{tabular}{|l|p{6cm}|p{4cm}|}
\hline
\textbf{Measure} & \textbf{Definition} & \textbf{Formula} \\
\hline
Support(\(X\)) & The empirical probability, \( P(X) \), of the itemset \( X \) & \( \frac{\text{N(X)}}{\text{Total number of transactions}} \) \\
\hline
Support(\(X \Rightarrow Y\)) & Proportion of transactions that contain both \(X\) and \(Y\) & \( \mathrm{supp}(X \cup Y) \) \\
\hline
Confidence(\(X \Rightarrow Y\)) & The conditional probability \( P(Y \mid X) \) & \(  \frac{\mathrm{supp}(X \cup Y)}{\mathrm{supp}(X)} \) \\
\hline
Lift(\(X \Rightarrow Y\)) & This measure evaluates the strength of the rule \(X \Rightarrow Y\)& \(  \frac{\mathrm{conf}(X \Rightarrow Y)}{\mathrm{supp}(Y)} \) \\
\hline
Conviction(\(X \Rightarrow Y\)) & This measure evaluates the violation of the rule \(X \Rightarrow Y\)  & \(  \frac{1 - \mathrm{supp}(Y)}{1 - \mathrm{conf}(X \Rightarrow Y)} \) \\
\hline
Leverage(\(X \Rightarrow Y\)) & It assesses the added value of the rule \(X \Rightarrow Y\) & \( \mathrm{supp}(X \cup Y) - \mathrm{supp}(X) \cdot \mathrm{supp}(Y) \) \\
\hline
\end{tabular}
\caption{Definitions and formulas of common statistical measures for association rules.}
\label{tab:assoc_metrics}
\end{table}

where \( N(X) = \text{Number of transactions containing } X\).
Hence, \[
P(X) = \frac{\text{Number of transactions containing } X}{\text{Total number of transactions}}.
\]

To apply association rules, we use the \textbf{Apriori algorithm}, a widely used technique in association rule mining. The Apriori method follows an iterative approach, starting with the identification of \textit{frequent itemsets}. 

\begin{enumerate}
    \item First, a \textbf{minimum support threshold} is defined.
    \item All itemsets with \textbf{support} greater than or equal to this threshold are identified. These are known as \textit{frequent itemsets} or \textit{large itemsets}.
    \item The algorithm then iteratively generates larger itemsets from smaller frequent itemsets, using the Apriori property which states that \textit{all subsets of a frequent itemset must also be frequent}.
    \item Once all frequent itemsets are identified, the next step is to generate \textbf{association rules} from them.
    \item For each rule of the form $X \Rightarrow Y$, where $X$ and $Y$ are itemsets and $X \cap Y= \emptyset$, the \textbf{confidence} is calculated as:
    \[
        \text{Confidence}(X\Rightarrow Y) = \frac{\text{Support}(X \cup Y)}{\text{Support}(X)}
    \]
    \item Only rules with confidence greater than or equal to the \textbf{minimum confidence threshold} are considered as \textbf{strong association rules}.
\end{enumerate}
These \textbf{strong association rules} are then considered as \textbf{the association rules of the data} under consideration. 

\section{Motivation}
\label{sec:motivation}

The computation of topological invariants  of Complete Intersection Calabi-Yau (CICY) manifolds is an active and important area of research. While the  invariants of CICY threefolds were computed in \cite{pgreen}, yielding 7890 distinct configuration matrices and a list of $921,497$ configuration matrices for CICY four-folds were worked out in ref. \cite{gray},  much less is known about their higher-dimensional analogues, such as Calabi-Yau 5-folds and 6-folds. These higher-dimensional CICYs are particularly interesting in M-theory and F-theory models, where compactification spaces of dimension five or more naturally arise. The size and complexity of configuration matrices grow rapidly, making  computations expensive using traditional algebraic geometry algorithms. 

The topological invariants of CICYs are not merely mathematical curiosities—they have direct physical implications. For instance, the Hodge numbers determine the number of moduli fields in the effective theory, while intersection numbers influence Yukawa couplings, anomaly cancellation conditions, and gauge symmetry breaking. Moreover, the so-called mirror symmetry provides dual pairs of Calabi-Yau manifolds whose Hodge diamonds exhibit symmetric structures. Understanding these symmetries in higher-dimensional settings remains an important open problem, both mathematically and physically.

Given the computational challenges of traditional methods, alternative such as data-driven approaches offer a powerful and scalable direction. In particular, \emph{association rule mining}—a machine learning technique can be adapted to uncover frequently co-occurring patterns among topological invariants. Instead of computing each invariant individually, this method searches for statistically significant rules of the form ``if $X$, then $Y$,'' where $X$ and $Y$ represent combinations of properties (e.g., Hodge number values or configuration features). Applying these techniques to higher-dimensional CICYs could provide new insights into their topological invariants, guiding further theoretical developments where conventional tools are no longer feasible.

\section{Results of Association Rules for CICY5 and CICY6 Datasets}
\label{sec:results}

To apply association rule mining to the datasets of complete intersection Calabi–Yau five-folds (\textbf{CICY5}) and six-folds (\textbf{CICY6}), we first transform the given Hodge numbers into a set of discrete items. The correspondence between the Hodge numbers and item labels is given in Table~\ref{tab:cicy-items}.

\begin{table}[H]
\centering
\renewcommand{\arraystretch}{1.5}
\begin{tabular}{|>{\centering\arraybackslash}m{2.5cm}|>{\centering\arraybackslash}m{3.5cm}||>{\centering\arraybackslash}m{2.5cm}|>{\centering\arraybackslash}m{3.5cm}|}
\hline
\multicolumn{2}{|c||}{\textbf{CICY5}} & \multicolumn{2}{c|}{\textbf{CICY6}} \\
\hline
\textbf{Item} & \textbf{Value} & \textbf{Item} & \textbf{Value} \\
\hline
item\(_1\) & \( h^{1,1} \) & item\(_1\) & \( h^{1,1} \) \\
item\(_2\) & \( h^{2,1} \) & item\(_2\) & \( h^{1,2} \) \\
item\(_3\) & \( h^{1,3} \) & item\(_3\) & \( h^{1,3} \) \\
item\(_4\) & \( h^{1,4} \) & item\(_4\) & \( h^{1,4} \) \\
item\(_5\) & \( h^{2,2} \) & item\(_5\) & \( h^{1,5} \) \\
item\(_6\) & \( h^{2,3} \) & item\(_6\) & \( h^{2,2} \) \\
--        & --             & item\(_7\) & \( h^{2,3} \) \\
--        & --             & item\(_8\) & \( h^{2,4} \) \\
--        & --             & item\(_9\) & \( h^{3,3} \) \\
\hline
\end{tabular}
\caption{Item encoding of Hodge numbers for CICY5 and CICY6 datasets.}
\label{tab:cicy-items}
\end{table}

The Apriori algorithm is applied to the dataset of complete intersection Calabi–Yau 6-folds (CICY6). The minimum support is set to \(10\%\), and the minimum confidence is set to \(80\%\). These threshold values are chosen based on the relative size of the dataset under consideration, ensuring a balance between rule significance and result comprehensiveness. The algorithm is implemented in \texttt{RStudio} using the \texttt{arulesViz} package and executed on the CCMS computing cluster.

The first ten results of the association rule mining applied to the CICY5 dataset are presented in Table~\ref{tab:association_rules5}, and there are 60 such rules in total. The remaining rules are provided in reference~\cite{Kaniba3}.

\begin{table}[H]
\centering
\caption{Top 10 Association Rules for CICY5}
\label{tab:association_rules5}
\begin{tabular}{|c|l|l|r|r|r|r|r|}
\hline
\textbf{rule} & \textbf{LHS} & \textbf{RHS} & \textbf{Support} & \textbf{Confidence} & \textbf{Coverage} & \textbf{Lift} & \textbf{Count} \\
\hline
1 & \{\} & \{item3=0\} & 0.9570 & 0.9570 & 1.0000 & 1.0000 & 11899 \\
2 & \{\} & \{item2=0\} & 0.9770 & 0.9770 & 1.0000 & 1.0000 & 12147 \\
3 & \{item5=4\} & \{item1=3\} & 0.1055 & 0.9661 & 0.1092 & 2.4962 & 1312 \\
4 & \{item5=4\} & \{item3=0\} & 0.1075 & 0.9838 & 0.1092 & 1.0280 & 1336 \\
5 & \{item5=4\} & \{item2=0\} & 0.1091 & 0.9993 & 0.1092 & 1.0228 & 1357 \\
6 & \{item5=7\} & \{item1=4\} & 0.1130 & 0.9256 & 0.1221 & 2.1287 & 1405 \\
7 & \{item5=7\} & \{item3=0\} & 0.1194 & 0.9776 & 0.1221 & 1.0215 & 1484 \\
8 & \{item5=7\} & \{item2=0\} & 0.1219 & 0.9987 & 0.1221 & 1.0222 & 1516 \\
9 & \{item5=5\} & \{item1=3\} & 0.1953 & 0.9870 & 0.1979 & 2.5501 & 2428 \\
10 & \{item5=5\} & \{item3=0\} & 0.1931 & 0.9760 & 0.1979 & 1.0198 & 2401 \\
\hline
\end{tabular}
\end{table}

In the same manner, the first ten results of the association rule mining applied to the CICY6 dataset are presented in Table~\ref{tab:association_rules6}. Our investigation yielded 160 association rules, the complete details of which are also provided in reference~\cite{Kaniba3}. 

\begin{table}[H]
\centering
\caption{Top 10 Association Rules for CICY6}
\label{tab:association_rules6}
\begin{tabular}{|c|l|l|r|r|r|r|r|}
\hline
\textbf{rule} & \textbf{LHS} & \textbf{RHS} & \textbf{Support} & \textbf{Confidence} & \textbf{Coverage} & \textbf{Lift} & \textbf{Count} \\
\hline
1  & \{\}           & \{item2=0\}         & 1.0000000 & 1.0000 & 1.0000000 & 1.0000 & 1482022 \\
2  & \{\}           & \{item3=0\}         & 1.0000000 & 1.0000 & 1.0000000 & 1.0000 & 1482022 \\
3  & \{\}           & \{item4=0\}         & 1.0000000 & 1.0000 & 1.0000000 & 1.0000 & 1482022 \\
4  & \{\}           & \{item7=0\}         & 1.0000000 & 1.0000 & 1.0000000 & 1.0000 & 1482022 \\
5  & \{item5=5\}    & \{item2=0\}         & 0.1319083 & 1.0000 & 0.1319083 & 1.0000 & 195491 \\
6  & \{item5=5\}    & \{item3=0\}         & 0.1319083 & 1.0000 & 0.1319083 & 1.0000 & 195491 \\
7  & \{item5=5\}    & \{item4=0\}         & 0.1319083 & 1.0000 & 0.1319083 & 1.0000 & 195491 \\
8  & \{item5=5\}    & \{item7=0\}         & 0.1319083 & 1.0000 & 0.1319083 & 1.0000 & 195491 \\
9  & \{item5=4\}    & \{item2=0\}         & 0.1691729 & 1.0000 & 0.1691729 & 1.0000 & 250718 \\
10 & \{item5=4\}    & \{item3=0\}         & 0.1691729 & 1.0000 & 0.1691729 & 1.0000 & 250718 \\
\hline
\end{tabular}
\end{table}

These results provide strong evidence of the usefulness of applying association rule mining to Calabi–Yau datasets, as they uncover statistically significant and interpretable relationships between different Hodge numbers.

\section{Conclusions and Future Outlooks}
\label{sec:conclusion}
In this paper, association rule machine learning is applied to datasets of complete intersection Calabi–Yau (CICY) 5-folds and 6-folds. For the 5-folds, we discover 60 significant association rules, while for the 6-folds, we identify 160 association rules across a dataset containing $1,482,022$ examples. The detailed rules are provided in Ref.~\cite{Kaniba3}. Our investigation reveals rich and non-trivial correlations among the Hodge numbers of these higher-dimensional Calabi–Yau manifolds. In particular, the  structural patterns (rules)  may guide future computations. Furthermore, this work demonstrates the potential of data-driven and interpretable machine learning methods in uncovering hidden relationships within topological invariants of CICYs datasets.
\section{Acknowledgment}
The author is indebted to Bobby Samir Acharya for pointing out the importance of CICY 6-folds as the next non-trivial even-dimensional class of complete intersection Calabi–Yau manifolds.

\end{document}